\documentclass[12pt,a4paper,]{article}

\usepackage{amssymb}
\usepackage{amsmath}
\usepackage{amsthm}
\usepackage{xypic}

\theoremstyle{plain}
\newtheorem{theorem}{Theorem}
\newtheorem{proposition}[theorem]{Proposition}

\theoremstyle{remark}

\theoremstyle{definition}
\newtheorem{definition}[theorem]{Definition}

\def\varinjlim_#1{\lim\limits_{\longrightarrow\atop{#1}}}

\def\Hom{\mathop{\rm Hom}\nolimits}
\def\id{\mathop{\rm id}\nolimits}

\def\diag{\mathop{\rm diag}\nolimits}
\def\pt{\mathop{\rm pt}\nolimits}

\begin{document}
\author{A. V. Ershov}
\title{Supplement\\ to the paper ''Floating bundles and their
applications''}
\date{}
\maketitle
This paper is the supplement to the section 2
of the paper "Floating bundles and their applications" \cite{pap1}.
Below we construct the denumerable set of extensions
of the formal group of geometric cobordisms
$F(x\otimes 1,1\otimes x)$
by the Hopf algebra $H=\Omega_U^*(Gr).$

Let $F_i(x\otimes 1,1\otimes x),\; i=1,2$
be formal groups over ring $R.$
Recall the following definition.
\begin{definition}
A homomorphism of formal groups
$\varphi \colon F_1\rightarrow F_2$
is a formal series $\varphi(x)\in R[[x]]$ such that
$\varphi(F_1(x\otimes 1,1\otimes x))=F_2(\varphi(x)\otimes 1,
1\otimes \varphi(x)).$
\end{definition}
Let $H$ be a Hopf algebra over ring $R$ with
diagonal $\Delta;$
let ${\frak F}_i(x\otimes 1,1\otimes x),\; i=1,2$ be
formal groups over $H.$
\begin{definition}
A homomorphism of formal groups over Hopf algebra $H$
$\Phi \colon {\frak F}_1\rightarrow {\frak F}_2$
is a formal series $\Phi(x)\in H[[x]]$ such that
$(\Delta \Phi)({\frak F}_1(x\otimes 1,1\otimes x))={\frak F}_2(\Phi(x)
\otimes 1,1\otimes\Phi(x)).$
\end{definition}
Note that $\varepsilon(\Phi)\colon
(\varepsilon \otimes \varepsilon)({\frak F}_1)
\rightarrow (\varepsilon \otimes \varepsilon)({\frak F}_2)$
is the homomorphism of the formal groups over the ring $R$
(where $\varepsilon$ is the counit of the Hopf algebra $H$).
We say that the homomorphism $\Phi$ covers the homomorphism
$\varepsilon(\Phi).$

Let $R$ be the ring $\Omega_U^*(\pt);$ let
$F(x\otimes 1,1\otimes x)\in R[[
x\otimes 1,1\otimes x]]$
be the formal group of geometric cobordisms.
Let $H$ be the Hopf algebra $\Omega_U^*(Gr).$
By definition, put
$\varphi^{(1)}(x)=x,\; \varphi^{(-1)}(x)=\theta(x)$
and $\varphi^{(n)}(x)=F(x,\varphi^{(n-1)}(x)),$
where $\theta(x)\in R[[x]]$ is the inverse element in $F.$
Clearly, that $\varphi^{(n)}\colon F\rightarrow F$
is the homomorphism for every $n\in \mathbb{Z}.$
Power systems were considered by S. P. Novikov
and V. M. Buchstaber in \cite{novbuch}.

Below for any $n\in \mathbb{Z}$ we construct the
extension ${\frak F}^{(n)}(x\otimes 1,1\otimes x)$
of $F(x\otimes 1,1\otimes x)$ by $H$
and the homomorphism $\Phi^{(n)}\colon {\frak F}\rightarrow
{\frak F}^{(n)}$ such that
\begin{itemize}
\item[(i)] ${\frak F}^{(1)}={\frak F};$
\item[(ii)] $\varepsilon(\Phi^{(n)})=\varphi^{(n)}.$
\end{itemize}

Let $X$ be a finite $CW$-complex.
Recall that the set of FBSP over $X$ is the semigroup
with respect to the multiplication of FBSP \cite{pap1}.
Let $n$ be a positive integer.
Let us take the product of the FBSP $\widetilde{Gr}_{k,kl}$
(over $Gr_{k,kl}$) with itself $n$ times. It is the FBSP
over $Gr_{k,kl}$ with a fiber $\mathbb{C}P^{k^n-1}\times
\mathbb{C}P^{l^n-1}.$ By $\widetilde{Gr}_{k,kl}^{(n)}$
denote the obtained FBSP.
Let $\widehat{Gr}_{k,kl}^{(n)}$ be the corresponding
bundle over $Gr_{k,kl}$ with fiber $\mathbb{C}P^{k^n-1}.$
Let $\widehat{Gr}^{(n)}=\varinjlim_{(k,l)=1}\widehat{Gr}_{k,kl}^{(n)}.$
We have the evident fiber maps $\widehat{Gr}_{k,kl}\rightarrow
\widehat{Gr}_{k,kl}^{(n)},\quad \lambda^{(n)} \colon \widehat{Gr}\rightarrow
\widehat{Gr}^{(n)}$ and
the following commutative diagrams ($(km,ln)=1$):
\begin{equation}
\begin{array}{ccc}
\widehat{Gr}_{km,klmn} & \rightarrow & \widehat{Gr}_{km,klmn}^{(n)} \\
\uparrow && \uparrow \\
\widehat{Gr}_{k,kl}\times \widehat{Gr}_{m,mn} &
\rightarrow & \widehat{Gr}_{k,kl}^{(n)}\times
\widehat{Gr}_{m,mn}^{(n)},\\
\end{array}
\end{equation}
\begin{equation}
\label{2}
\begin{array}{ccc}
\widehat{Gr} & \stackrel{\lambda^{(n)}}{\rightarrow} & \widehat{Gr}^{(n)}\\
\scriptstyle{\widehat{\phi}}\uparrow && \uparrow
\scriptstyle{\widehat{\phi}^{(n)}} \\
\widehat{Gr}\times \widehat{Gr} & \stackrel{\lambda^{(n)}
\times \lambda^{(n)}}{\rightarrow} &
\widehat{Gr}^{(n)}\times \widehat{Gr}^{(n)}.\\
\end{array}
\end{equation}
By $x$ denote the class of cobordisms in
$\Omega_U^{2}(\widehat{Gr}^{(n)})$
such that its restriction to any fiber $\cong \mathbb{C}P^\infty$
is the generator
$x\mid _{\mathbb{C}P^\infty}\in \Omega_U^2(\mathbb{C}P^\infty)$.
Let $\Phi^{(n)}(x)\in H[[x]]$ be the series, defined by
the fiber map $\lambda^{(n)}$ (see \cite{pap1}).
Let
$$
{\frak F}^{(n)}(x\otimes 1,1\otimes x)\in H{\mathop{\widehat{\otimes}}
\limits_R}H[[x\otimes 1,1\otimes x]]
$$
be the series,
corresponds to the fiber map $\widehat{Gr}^{(n)}\times
\widehat{Gr}^{(n)}\stackrel{\widehat{\phi}^{(n)}}
{\rightarrow}\widehat{Gr}^{(n)}$
(see \cite{pap1}; note that $Gr^{(n)}$ is the $H$-group
with the multiplication $\widehat{\phi}^{(n)}$).
Clearly, that ${\frak F}^{(n)}(x\otimes 1,1\otimes x)$
is an extension of $F(x\otimes 1,1\otimes x)$ by $H$
(in particular, it is the formal group over Hopf algebra $H$).
Note that $\lambda^{(n)}$ covers the identity map
of the base $Gr.$ It follows from diagram (\ref{2}) that
$$
(\Delta \Phi^{(n)})({\frak F}(x\otimes 1,1\otimes x))=
{\frak F}^{(n)}(\Phi^{(n)}(x)\otimes 1,1\otimes \Phi^{(n)}(x)).
$$
It is clear that $\varepsilon(\Phi^{(n)})(x)=\varphi^{(n)}(x).$

For $n=0$ let $\widehat{Gr}^{(0)}=Gr\times \mathbb{C}P^{\infty}$
and let $\lambda^{(0)}$ be the composition
$$\widehat{Gr}\rightarrow \pt \rightarrow \widehat{Gr}^{(0)}.$$
It defines the series ${\frak F}^{(0)}=F$ and $\Phi^{(0)}=0.$

Let $\lambda^{(-1)}$ be the fiber map $\widehat{Gr}\rightarrow
\widehat{Gr}^{(-1)}=\widehat{Gr}$ such that the following conditions
hold:
\begin{itemize}
\item[(i)] the restriction of $\lambda^{(-1)}$ to any fiber
is the inversion in the $H$-group $\mathbb{C}P^{\infty}$
(i. e. the complex conjugation);
\item[(ii)] $\lambda^{(-1)}$ covers the map $\nu \colon Gr \rightarrow Gr,$
where $\nu$ is the inversion in the $H$-group $Gr.$
\end{itemize}
Let $\Phi^{(-1)}(x)\in H[[x]]$ be the series,
defined by $\lambda^{(-1)}.$
Trivially, that $\varepsilon(\Phi^{(-1)})(x)=\theta(x).$
Note that the $\lambda^{(-1)}$ coincides with $\widehat{\nu}$
(see \cite{pap1}). Consequently, $\Phi^{(-1)}=\Theta(x).$
Now we can define ${\frak F}^{(n)}$ and $\Phi^{(n)}$
for negative integer $n$ by the obvious way.

By $S$ denote the antipode
of the Hopf algebra $H.$ Let $\mu$ be the multiplication
in the Hopf algebra $H.$
By definition, put $(1)=\id_H,\; (-1)=S\colon H\rightarrow H$ and
$(n)=\mu \circ ((n-1)\otimes (1))\circ \Delta \colon H\rightarrow H$
(in particular, $(0)=\eta \circ \varepsilon \colon H\rightarrow H,$
where $\eta$ is the unit in $H$).
\begin{proposition}
${\frak F}^{(n)}(x\otimes 1,1\otimes x)=
(((n)\otimes(n)){\frak F})(x\otimes 1,1\otimes x)$
for any $n\in \mathbb{Z}.$
\end{proposition}
{\raggedright {\it Proof}.}
By $\phi \colon Gr\times Gr\rightarrow Gr$ denote the multiplication
in the $H$-space $Gr.$
Suppose $n$ a positive integer.
By definition, put $\phi{(1)}=\id_{Gr},\;
\phi{(n)}=\phi \circ (\phi{(n-1)} \times \id_{Gr}),$
and $\diag{(n)}=(\diag{(n-1)}\times \id_{Gr})\circ \diag ,$
where $\diag(1)=\id_{Gr},\;
\diag =\diag{(2)}\colon Gr\rightarrow Gr\times Gr.$
Note that the composition $\phi{(n)}\circ \diag{(n)}\colon Gr
\rightarrow Gr$ induces the homomorphism $(n)\colon H\rightarrow H.$

Let us consider the classifying
map $\alpha{(n)}\colon Gr\rightarrow Gr$ for the bundle $\widehat{Gr}^{(n)}$
over $Gr.$ We have the following commutative diagram:
$$
\diagram
\mathbb{C}P^\infty \drto \rto^= & \mathbb{C}P^\infty \drto \\
& \quad \widehat{Gr}^{(n)} \dto \rto^{\widehat{\alpha}(n)} &
\quad \widehat{Gr} \dto \\
& Gr \rto^{\alpha{(n)}} & Gr \\
\enddiagram
$$
It is easy to prove that $\alpha{(n)}=\phi{(n)}\circ \diag{(n)}.$
Hence $\alpha{(n)}^*=(n)\colon H\rightarrow H.$
Note that the following diagram
\begin{equation}
\begin{array}{ccc}
\widehat{Gr}^{(n)}\times \widehat{Gr}^{(n)} &
\stackrel{\widehat{\alpha}{(n)}
\times \widehat{\alpha}{(n)}}{\rightarrow} &
\widehat{Gr}\times \widehat{Gr}\\
\scriptstyle{\widehat{\phi}^{(n)}}\downarrow \qquad &&
\quad \downarrow \scriptstyle{\widehat{\phi}} \\
\widehat{Gr}^{(n)} & \stackrel{\widehat{\alpha}{(n)}}
{\rightarrow} & \widehat{Gr}\\
\end{array}
\end{equation}
is commutative. This completes the proof for positive $n.$
For negative $n$ proof is similar. $\square$

\smallskip

We can define the structure of group on the set $\{ {\frak
F}^{(n)};\; n\in \mathbb{Z}\}$ in the following way.
Recall that for any Hopf algebra $H$ the triple $(\Hom_{Alg. Hopf}(H,H),
\star, \eta \circ \varepsilon)$
is the algebra with respect to the convolution $f\star g=
\mu \circ (f\otimes g)\circ \Delta \colon H\rightarrow H.$
It follows from the previous Proposition that the formal group
${\frak F}^{(n)}$ corresponds to the homomorphism $(n)\colon H
\rightarrow H$ (see Conjecture 24 in \cite{pap1}).
Clearly, that $(m)\star(n)=(m+n)$ for any $m,n\in \mathbb{Z}.$

\end{document}